\input amstex   
   
\documentstyle{amsppt}   
\mag=\magstep1   
\NoRunningHeads   
\NoBlackBoxes   
\topmatter   
\title Some curvature estimates for Riemannian manifolds\\ equipped with  
foliations of rank 2\endtitle   
\author Pawe\l \/ Walczak\endauthor   
\address Instytut Matematyki, Uniwersytet \L\'odzki, Banacha 22, 90238    
\L\'od\'z, Poland\endaddress   
\email  pawelwal\@ krysia.uni.lodz.pl\endemail   
\date  December 30, 1995\enddate   
\thanks Supported by the KBN grant no 2P30103604\endthanks
\subjclass 53 C 20\endsubjclass
\abstract Some curvature estimates are derived from geometrical data 
concerning quasi-conformality properties of some commuting linearly independent
vector fields on a compact Riemannian manifold.
\endabstract
\endtopmatter   
   
\document   
   
\subhead 1. Introduction\endsubhead   
Calculating the rank, i.e. the maximal number of pointwise linearly independent    
commuting vector fields, of a manifold $M$ is a well known problem in differential    
topology. If rank$\, M = k$, then $M$ admits a locally free action of $\Bbb R^k$.  
Such   actions (and foliations generated by them) are classified up to the either  
topological or   differentiable type in several cases ([{\bf AC}], [{\bf AS}], [{\bf  
CRW}], [{\bf RR}],   etc.). Also, for some pairs $(k, n)$, $n$-dimensional closed  
manifolds of rank $k$ are   classified. For example, if $k = n-1$, then $M$ is a {\bf  
T}$^m$-bundle over a torus   {\bf T}$^{n-m}$ for some $m$ ([{\bf CR}], [{\bf R}],  
[{\bf RRW}], etc.).   
   
On the other hand, in several geometric situations (like infinitesimal    
isometries [{\bf Be}] or infinitesimal conformal transformations on closed Riemannian  
manifolds of positive sectional curvature, see [{\bf W}] and    
[{\bf Y]}, for instance) vector fields have to    
vanish at some points or have to be somewhere collinear whenever commute. The result    
presented here is of this sort: Assuming the existence of two commuting linearly    
independent vector fields $X$ and $Y$ on a Riemannian manifold $(M,g)$ we    
estimate the curvature of $M$ in terms of some geometric data (like the divergence, the    
norm of  the Ahlfors operator) obtained from the $\Bbb R^2$-action generated by the flows of $X$ and $Y$. More precisely, our result    
reads as follows.   
   
\proclaim{Theorem} For any $\epsilon > 0$  there exists a constant    
$K(\epsilon) > 0$ such that any closed $n$-dimensional ($n\ge 2$) 
Riemannian manifold of sectional curvature $K_M > K(\epsilon )$ 
admits no $(C^1,\epsilon )$-quasi-conformal $\Bbb R^2$-actions.
\endproclaim   
   
Roughly speaking, a vector field $X$ is $C^1$-quasi-conformal if the Ahlfors
operator $S$ on $(M, g)$ evaluated at $X$ has the $C^1$-norm small enough. 
An action $\phi$ of $\Bbb R^2$ is $(C^1,\epsilon)$-quasi-conformal if all
the vector fields $Z_u$ ($u\in\Bbb R^2$) given by $Z_u(p) =
(t\mapsto \phi (tu, p)){\dot{}} (0)$ ($p\in M$) are $(C^1,\epsilon )$-quasi-conformal. For the precise definition of  infinitesimal $(C^1, \epsilon )$-quasi-conformal transformations, see Section 3.

The proof of the Theorem is given in Sections 2 - 5 below. In general, 
it follows the line  of that in [{\bf Y}]. Several 
identities in there which hold for infinitesimal conformal transformations 
are replaced by inequalities produced from the geometric    
data related to the vector fields under consideration. At some points,
the argument differs essentially. For instance, instead of studying the
eigenvalues and eigenvectors of some linear operator $A$ we estimate the inner
products of $A(v)$ with some particular vector $y$. Section 6 contains an application 
of the estimates derived before: We obtain an inequality which involves 
the Gaussian curvature of the orbits of a locally free  
$(C^1, \epsilon )$-quasi-conformal action of $\Bbb R^2$ and their second fundmental
form.
   
\demo{Remark} Our Theorem together with the mentioned above results on the rank of closed  manifolds justify to some extend the following Yang's remark [{\bf Y}]: 
One might expect that any two commuting vector fields on a closed manifold 
of positive sectional curvature must be linearly dependent at some points, 
in other words, if $M$ is closed and has positive sectional curvature, then 
rank$\, M\le 1$. However, the problem of calculating (or, estimating) the rank 
of all manifolds which admit positively curved Riemannian metrics is probably
pretty difficult. First, there is a lot of examples of positively curved closed 
manifolds of different topological types ([{\bf AW}], [{\bf E1}], [{\bf E2}], 
[{\bf KS}], etc.) and still new examples arrive (see [{\bf Ba}] and [{\bf T}]
for instance). Second, even if $M$ is simply connected, the curvature of $M$ is 
$\delta$-pinched with $\delta$ close enough to 1 (for example, $\delta = 0.681$) 
and, therefore, $M$ is diffeomorphic to the standard sphere $S^n$ [{\bf S}], 
the rank of $M$ is unknown in odd dimensions. Of course, rank$\, S^n = 0$ when 
$n$ is even. Moreover, rank$\, S^3 = 1$ [{\bf L}] and there are no locally free 
actions of $S^1\times\Bbb R$ on any standard sphere [{\bf ASc}].  
\enddemo

\subhead 2. A preliminary lemma\endsubhead
Let $X$ be an arbitrary vector field on a Riemannian manifold $(M, g)$, 
$p_0\in M$ and $X(p_0)\ne 0$.

\proclaim{Lemma} For any vector $v\in T_{p_0}M$, $v\ne 0$, orthogonal to 
$X(p_0)$ there exists a vector field $V$ in a neighbourhood $U$ of $p_0$  
such that 
$$V(p_0) = v,\ \langle V, X\rangle = 0,\ [V, X] = \theta\cdot X\ \ \text{and}\ \ 
\nabla_vV = a\cdot X(p_0)\tag1$$
for some $\theta\in C^{\infty}(U)$ and $a\in\Bbb R$.
\endproclaim

\demo{Proof} In a neighbourhood $U$ of $p_0$, there exists a chart 
$\phi = (x, y, z_1,\dots , z_{n-2})$ ($n = \dim M$) such that 
$X = \frac{\partial}{\partial x}$, $v = \frac{\partial}{\partial y}(p_0)$ and 
$\nabla_v\frac{\partial}{\partial y} = 0$. Find a function $f\in C^{\infty}(U)$ such that 
$V = \frac{\partial}{\partial y} + f\cdot\frac{\partial}{\partial x}$  is orthogonal to 
$X$. $V$ satisfies (1) with
$$\theta = -\frac{\partial f}{\partial x}\ \ \text{and}\ \ 
a = \frac{\partial f}{\partial y}(p_0).\qed\tag2$$
\enddemo

\subhead 3. Infinitesimal $(C^1, \epsilon)$-quasi-conformal transformations\endsubhead   
Given $\epsilon\ge 0$, a vector field $X$ on a Riemannian manifold $(M, g)$ will be    
called an infinitesimal $(C^1, \epsilon)$-quasi-conformal transformation 
(an infinitesimal  $(C^1, \epsilon)$-QCT,  for short) if    
$$\| SX(p)\|\le\epsilon\| X(p)\| \ \ \text{and}\ \  \|\nabla SX(p)\|\le\epsilon^2\| 
X(p)\|\tag 3$$ 
at any point $p$ of $M$. 

Here $S$ is the Ahlfors operator on $M$.
Recall that $SX$ is a traceless symmetric 2-form on $M$ defined by    
$$SX = \Cal L_Xg - f_X\cdot g,\tag 4$$   
where $\Cal L_X$ denotes the Lie derivation with respect to $X$, 
$f_X = \frac2{n}\text{div}X$ and $n = \dim M$. Equivalently, $SX$ can be defined by   
$$SX(V,W) = \langle\nabla_VX, W\rangle + \langle V, \nabla_WX\rangle -    
f_X\cdot\langle V, W\rangle,\tag 5$$   
where $V$ and $W$ are arbitary vector fields and $\nabla$ is the Levi-Civita    
connection on $(M, g)$. Recall also that $\| SX\|$ describes the rate of  
quasi-conformality of the flow maps $(\phi_t)$ of $X$, namely, any map $\phi_t$ is  
$k_t$-quasi-conformal with $k_t\le \exp ({|t|\cdot\| SX\|})$, $\|\cdot\|$ being the  
supremum norm here [{\bf P}].  

\demo{Remark} Observe that our notion of quasi-conformality differs from the
standard one. Usually, a vector field $X$ is said to be $k$-quasi-conformal just
when $\| SX\|\le k$. Also, observe that if the $C^1$-norm of $\nabla X$ is small,
$$\|\nabla X\|\le\delta\cdot\| X\|\ \ \text{and}\ \ 
\|\nabla^2X\|\le\delta^2\cdot\| X\| ,\tag 6$$
then $X$ is $(C^1, \epsilon )$-quasi-conformal with $\epsilon = 4\delta$. 
Finally, note that if $X$ is an infinitesimal $(C^1, \epsilon)$-QCT then, 
for any $c\ne 0$, $cX$ is an infinitesimal $(C^1, \epsilon)$-QCT as well, and 
if $X$ is an infinitesimal $(C^1, \epsilon)$-QCT   
on $(M, g)$ and $g' = c^2g$ for some $c\in\Bbb R_{+}$, then $X$ is an 
infinitesimal $(C^1, \epsilon ')$-QCT on $(M, g')$ with $\epsilon ' = 
c^{-1}\epsilon$.
\enddemo

\subhead 4. First estimates  \endsubhead   
Assume now that $X$ and $Y$ are pointwise linearly independent commuting    
vector fields on a compact Riemannian manifold $M$ such that $u_1X + u_2Y$ is an infinitesimal $(C^1, \epsilon)$-QCT for any $u_1$ and $u_2\in\Bbb R$.
 
For any $p\in  M$ and $t\in S^1 = \Bbb  R/2\pi\Bbb Z$ put   
$$Z(p,t) = \cos t\cdot X(p) + \sin t\cdot Y(p).\tag7$$   
The function $h:M\times S^1\to\Bbb R$ given by $h(p, t) = \| Z(t)(p)\|^2$ is    
continuous, therefore it attains its minimum at a point $(p_0, t_0)$ of $M\times S^1$.    
Since all the fields $Z(\cdot , t)$ are $(C^1, \epsilon )$-quasi-conformal,
we may assume without losing generality that that $t_0 = 0$. 

Let $x = X(p_0)$, $y = Y(p_0)$ and $v\in T_{p_0}M$ be an unit vector. The    
standard analysis of partial derivatives of $h$ at $p_0$ shows that   
$$v\langle X, X\rangle  = 0,\quad \langle x, y\rangle  = 0,\tag 8$$   
$$V^2\langle X, X\rangle (p_0)\ge 0,\tag9$$   
$$\| y\|\ge \| x\| \tag10$$   
and   
$$d = \vmatrix V^2\langle X, X\rangle (p_0) & v\langle X, Y\rangle \\   
v\langle X, Y\rangle &  \| y\|^2 - \| x\|^2\endvmatrix\ge 0,\tag11$$   
where $V$ is an arbitrary vector field on a neighbourhood of $p_0$, $V(p_0) = v$.   

Now, take any unit vector $v\perp x$ and extend it to a vector field $V$ satisfying the conditions of Lemma in Section 2. Then, at the point $p_0$ as above, 
$$SX(x, x) = - f_X(p_0)\cdot\| x\|^2,\ \ \text{therefore}\ \ |f_X(p_0)|\le\epsilon\| x\| ,\tag12$$
$$SX(x,v) =  - x\langle X,V\rangle - \langle X, [X,V]\rangle = \theta (p_0)\| x\|^2,\ \ \text{therefore}\ \ |\theta (p_0)|\le\epsilon\tag13$$
and
$$SX(v,v) = -2\langle x,\nabla_vV\rangle - f_X(p_0) = -2a\| x\|^2 -f_X(p_0),\ \ \text{therefore}\ \ a\| x\|\le\epsilon .\tag14$$

Moreover, $\langle\nabla_xX, x\rangle = 0$ and $\langle\nabla_xX,v\rangle = - \langle x, \nabla_xV\rangle = -\langle x,\nabla_vX\rangle - \langle x, [X,V](p_0)\rangle $ $= \theta(p_0)\cdot\| x\|^2$. From (13) it follows that
$$\|\nabla_xX\|\le\epsilon\| x\|^2 .\tag15$$

Also, since $SX(v,v) = x\langle V, V\rangle - f_X(p_0)$, we get from (3) and (12) that
$$|X\langle V,V\rangle (p_0)|\le2\epsilon\| x\| .\tag16$$

Passing to the covariant derivative $\nabla SX$ we observe that
$$\aligned (\nabla_xSX)(v,v) &= X^2\langle V, V\rangle (p_0) -
Xf_X(p_0) - f_X(p_0)X\langle V,V\rangle (p_0)
\\ &- 2SX(\nabla_xV,v),\endaligned\tag17$$
$$(\nabla_vSX)(x,v) = v\theta\cdot\| x\|^2 - SX(\nabla_xV,v) - \theta (p_0)^2\| x\|^2 - a\cdot SX(x,x),\tag18$$
and
$$(\nabla_xSX)(x,x) = X^2\langle X, X\rangle (p_0) - Xf_X(p_0)\cdot \| x\|^2 - 2SX(\nabla_xX, x).\tag19$$
Combining identities (17) - (19) and using inequalitites (3) and (12) - (16) we arrive at
$$\left| X^2\langle V,V\rangle (p_0) - 2v\theta\cdot\| x\|^2 - \| x\|^{-2}X^2\langle X,X\rangle (p_0)\right|\le 12\epsilon^2\| x\|^2.\tag20$$

All the inequalities above will be applied in the next Sections to get some 
curvature estimates.

\subhead 5. Estimates involving sectional curvature\endsubhead   
Let us keep all the notation introduced in Section 4 and assume that the 
sectional curvature of $M$ is positive. (Otherwise, the statement of the 
Theorem holds trivially.) Also, denote by $A$ the linear transformation of $T_{p_0}M$ given by
$$Aw = \nabla_wX\tag21$$ 
and choose a vector $v$ in such a way that 
$$Av = by\ \ \text{for some}\ \ b\in\Bbb R.\tag22$$
(Such a vector exists. Since $A:\{ x\}^\perp\to\{ x\}^\perp$, either Ker $A|\{ x\}^\perp$ is nontrivial and $Av = 0$ for some $v\perp x$ or $A|\{ x\}^\perp$ is an isomorphism.)
   
As before, extend $v$ to a vector field $V$ defined in a neighbourhood of 
$p_0$ as in our Lemma. From the definition of the curvature   
tensor $R$ on  $(M, g)$ we obtain easily the identity   
$$\aligned V^2\langle X, X\rangle (p_0) + X^2\langle V, V\rangle (p_0) &
= -2 K_M(x\wedge v)\| x\|^2 + 2\| Av|\|^2 \\ & + 2v\theta\cdot\| x\|^2 
- 2\theta (p_0)\langle Ax, v\rangle,\endaligned\tag23$$    
where   
$$K_M(x\wedge v) = \| x\|^{-2}\langle R(v,x)x, v\rangle\tag24$$   
is the sectional curvature of $M$ in the direction of the plane spanned by 
$x$ and $v$. Since both $X$ and $Y$ are infinitesimal  
$(C^1, \epsilon )$-QCTs,   
$$|\langle\nabla_xY, v\rangle + \langle x, \nabla_vY\rangle |\le\epsilon\| x\|\cdot\|   
y\|\tag25$$   
and   
$$|\langle\nabla_yX, v\rangle + \langle y, \nabla_vX\rangle |\le 2\epsilon\| x\|\cdot \|   
y\|.\tag26$$   
Since $X$ and $Y$ commute, $\nabla_xY = \nabla_yX$ and (25) together with (26)    
imply that   
$$|v\langle X, Y\rangle - 2 b\| y\|^2|\le 3\epsilon\| x\|\cdot\| y\| .\tag27$$   
Consequently,  
$$- (v\langle X, Y\rangle )^2\le - 4b^2\| y\|^4 + 12\epsilon b\| x\|\cdot\| y\|^3 +   
9\epsilon^2\| x\|^2\| y\|^2 .\tag28$$  

The determinant $d$ in (11) satisfies   
$$\aligned d = &(|y|^2 - \| x\|^2)(-2K(x\wedge z)\| x\|^2  - X^2\langle V, 
V\rangle (p_0) + 2\| Av\|^2\\ & - 2\theta (p_0)\langle Ax, v\rangle + 2v\theta\cdot\| x\|^2) - (v\langle X, Y\rangle )^2.\endaligned\tag 29$$   
Since the sectional curvature of $M$ is positive, combining inequalities (9), 
(13), (15), (20) and (28), and using identity (22)
we obtain the inequality    
$$\aligned 0\le d\le &(\| y\|^2 - \| x\|^2)(14\epsilon^2\| x\|^2 
+ 2b^2\| y\|^2) - b^2\| y\|^4 \\ &+ 
12\epsilon b\| x\|\cdot\| y\|^3 + 9\epsilon^2\| x\|^2\| y\|^2.  
\endaligned\tag 30$$    
Substituting $s = b\| x\|^{-1}\| y\|$, dividing both sides of (30) 
by $\| y\|^2$ and performing some other elementary transformations one can 
see easily that (30) implies  
$$-2s^2 + 12\epsilon s + 23\epsilon^2\ge 0.\tag31$$  
Therefore,  
$$s\le (3 + \frac12\sqrt{82})\epsilon\| x\|\tag32$$  
and 
$$b\le (3 + \frac12\sqrt{82})\epsilon\| x\|\cdot\| y\|^{-1} .\tag33$$  

Applying formula (22) once again we can see that  
$$\aligned K(x\wedge v) &  = \| x\|^{-2}( -\frac12V^2\langle X, X\rangle    
(p_0) - \frac12X^2\langle V, V\rangle (p_0) + \| Av\|^2 \\& + 
v\theta - \theta (p_0)\langle Ax, v\rangle )
\le \| x\|^{-2}\left( 6\epsilon^2\| x\|^2 + 
b^2\| y\|^2\right)\\ &\le\epsilon^2(36.5 + 3 \sqrt{82})
\approx 63.66615541\dots\epsilon^2\endaligned\tag34$$  
Therefore, the statement of the Theorem holds with    
$$K(\epsilon ) = \epsilon^2 (36.5 + 3\sqrt{82}).\qed\tag35$$  
 
\proclaim{Corollary [Y]} Any two commuting conformal vector fields on a closed  
Riemannian manifold of positive sectional curvature must be linearly dependent 
at least at one point. 
\endproclaim

\subhead 6. Estimates involving Gaussian curvature of the leaves\endsubhead Let $\Cal F$ be a foliation of a compact  
positively curved manifold $M$ defined by a locally free action of $\Bbb R^2$. 
Then $\dim F = 2$ and the tangent bundle of $\Cal F$ is spanned by two
commuting linearly  independent vector fields $X$ and $Y$. These fields as well
as all the linear combinations $u_1X + u_2Y$, $u_i\in\Bbb R$, become infinitesimal
$(C^1, \epsilon )$-QCTs for some $\epsilon\ge 0$. Therefore, with this $\epsilon$,
$\Cal F$ is defined  by a $(C^1,\epsilon )$-quasi-conformal action of $\Bbb R^2$. 
 
Consider fields $X$ and $Y$, and find $p_0\in M$ as before. 
Let $v$ be a unit vector tangent to $\Cal F$ and orthogonal to $x = X(p_0)$,
say $v = y/\| y\|$. Then 
$$Av = by + B(x, v),\tag36$$ 
where $B$ is the second fundamental form of $\Cal F$.   
The argument analogous to that  of Section 5 shows that 
$$b\le\| x\|\cdot\| y\|^{-1}\left((3 + \frac{1}{2}\sqrt{82})\epsilon +  
\| B (p_0)\|\right) ,\tag37$$ 
where $x_1 = x/\| x\|$. Applying formula (22) as in Section 5 we can get the estimate 
$$K_M(x\wedge v)\le\epsilon^2(36.5 + 3\sqrt{82}) + \epsilon\sqrt2 (3 
+ 0.5\sqrt{82})\|  B(p_0)\| + \frac32\| B(p_0)\|^2.\tag38$$ 
 
Finally, the classical Gauss equation reads
$$K_{\Cal F}(p_0) = K_M(x\wedge v) + \langle B(x_1, x_1), B(v,v)\rangle 
- \| B(x_1, v)\|^2,\tag39$$ 
where $K_{\Cal F}$ is the Gaussian curvature of the leaves of $\Cal F$
and $x_1 = x/\| x\|$. Therefore, we can conclude by the following. 
 
\proclaim{Proposition} The Gaussian curvature $K_{\Cal F}$ of a foliation $\Cal F$  
defined by a locally free $(C^1, \epsilon )$-quasi-conformal action of 
$\Bbb R^2$ on a closed positively curved Riemannian manifold $M$ satisfies 
at some point $p_0$ of $M$ the inequality 
$$K_{\Cal F}(p_0)\le\epsilon^2(36.5 + 3\sqrt{82}) + 
(6 + \sqrt{82})\|  B(p_0)\| + \| B(p_0)\|^2\tag40$$ 
with $B$ being the second fundamental form of $\Cal F$.\qed 
\endproclaim 
   
\enddocument